\documentclass[12pt,thmsa]{article}
\usepackage{amsmath, latexsym, amsfonts, amssymb, amsthm, amscd}

\newtheorem{theorem}{Theorem}
\newtheorem{corollary}{Corollary}
\newtheorem{proposition}{Proposition}
\newtheorem{lemma}{Lemma}

\newtheorem{defi}{Definition}
\newcommand{\p}{\Bbb{P}}

\newcommand{\e}{\Bbb{E}}

\newcommand{\ind}{\mbox{\rm 1\hspace{-0.04in}I}}

\newcommand{\R}{\mbox{\rm I\hspace{-0.02in}R}}

\newcommand{\ed}{\stackrel{(d)}{=}}
\newcommand{\eqdef}{\stackrel{\mbox{\tiny$($def$)$}}{=}}

\def\QED{\hfill\vrule height 1.5ex width 1.4ex depth -.1ex \vskip20pt}

\begin{document}

\title{On the genealogy of conditioned\\ stable L\'evy forests.}

\maketitle

\begin{center}
{\large L. Chaumont}\footnote{LAREMA, D\'epartement de
Math\'ematiques, Universit\'e d'Angers, 2, Bd Lavoisier - 49045,
{\sc Angers Cedex 01.}  E-mail: loic.chaumont@univ-angers.fr}{\large
 and  J.C. Pardo}\footnote{Laboratoire de Probabilit\'es et Mod\`eles
Al\'eatoires, Universit\'e Pierre et Marie Curie, 4, Place Jussieu -
75252 {\sc Paris Cedex 05.} E-mail: pardomil@ccr.jussieu.fr}
\end{center}
\vspace{0.2in}

\begin{abstract} We give a realization of the stable L\'evy forest
of a given  size conditioned by its mass from the path of the
unconditioned forest. Then, we prove an invariance principle for
this conditioned forest by considering $k$ independent Galton-Watson
 trees whose offspring distribution is in the domain of
attraction of any stable law conditioned on their total progeny to
be equal to $n$. We prove that  when $n$ and $k$ tend towards
$+\infty$, under suitable rescaling, the associated coding random
walk, the contour and height processes converge in law on the
Skorokhod space respectively towards the ``first passage bridge'' of
a stable L\'evy
process with no negative jumps and its height process.\\


\end{abstract}

\newpage

\section{Introduction}
\noindent The purpose of this work is to study some remarkable
properties of  stable L\'evy forests of a given size conditioned by their
 mass.\\

A Galton-Watson tree is the underlying family tree of a given
Galton-Watson process with offspring distribution $\mu$ started with
one ancestor.  It is well-known that if $\mu$ is critical or
subcritical, the Galton-Watson process is almost surely finite and
therefore, so is the corresponding Galton-Watson tree. In this case,
Galton-Watson trees can be coded by two different discrete real
valued processes: the height process and the contour process whose
definition is recalled here in section 2. Both processes describe
the genealogical structure of the associated Galton-Watson process.
They are not Markovian but can be written as functionals of a
certain left-continuous random walk whose jump distribution depends
on the offspring distribution $\mu$. In a natural way, Galton-Watson
forests are a finite or infinite
collections of independent Galton-Watson trees.\\

The definition of  L\'evy trees  bears upon on the continuous
analogue of the height process of Galton-Watson trees introduced by
Le Gall and Le Jan in \cite{LL}  as a functional of a L\'evy process
with no negative jumps. Our presentation owes a lot to the recent
paper of Duquesne and Le Gall \cite{DL2}, which uses the formalism
of $\R-$trees to define  L\'evy trees that were implicit in
\cite{Du}, \cite{DL} and \cite{LL}. We may consider L\'evy trees as
random variables taking values in the space of compact rooted
$\R$-trees. In a recent paper  of Evans, Pitman and Winter
\cite{EPW}, $\R-$trees are studied from the point of view of mesure
theory. Informally an $\R-$tree is a metric space $(\mathcal{T}, d)$
such that for any two points $\sigma$ and $\sigma'$ in $\mathcal{T}$
there is a unique arc with endpoints $\sigma$ and $\sigma'$ and
furthermore this arc is isometric to a compact interval of the real
line. In \cite{EPW}, the authors also established that the space
$\mathbb{T}$ of equivalent classes of (rooted) compact real trees,
endowed with the Gromov-Hausdorff metric, is a Polish space. This
makes it very natural to consider random variables or even random
processes taking values in the space $\mathbb{T}$. In this work, we
define L\'evy forests as Poisson point processes with values in the
set of $\R$-trees whose characteristic measure is the law of the
generic L\'evy tree.\\

First, we are interested in the construction of L\'evy forests of a
given size conditioned by their mass. Again, in the discrete setting
this conditioning is easier to define; the conditioned Galton-Watson
forest of size $k$ and mass $n$ is a collection of $k$ independent
Galton-Watson trees with total progeny equal to $n$. In section 4,
we provide a definition of these notions for L\'evy forest. Then, in
the stable case, we give a construction of the conditioned stable
L\'evy forest of size $s>0$ and  mass $1$ by rescaling the
unconditioned forest of a particular
random mass.\\

In \cite{Al}, Aldous showed that the Brownian random tree (or
continuum random tree) is the limit as $n$ increases of a rescaled
critical Galton-Watson tree conditioned to have $n$ vertices whose
offspring distribution has a finite variance. In particular, Aldous
proved that the discrete height process converges on the Skorokhod
space of c\`adl\`ag paths to the normalized Brownian excursion.
 Recently, Duquesne \cite{Du} extended such results
to Galton-Watson trees whose offspring distribution is in the domain
of attraction of a stable law with index $\alpha$ in $(1, 2]$. Then,
Duquesne showed that the discrete height process of the Galton
Watson tree conditioned to have a deterministic progeny, converges
as this progeny tends to infinity on the  Skorokhod space to the
normalized excursion of the height process associated with the
stable L\'evy process.\\

 The other main purpose of our work is to study
this convergence in the case of a finite number of independent
Galton-Watson trees, this number being an increasing function of the
progeny. More specifically, in Section 5, we establish an invariance
principle for the conditioned forest by considering $k$ independent
Galton-Watson trees whose offspring distribution is in the domain of
attraction of any stable law conditioned on their total progeny to
be equal to $n$. When $n$ and $k$ tend towards $\infty$, under
suitable rescaling, the associated coding random walk, the contour
and height processes converge in law on the space of Skorokhod
towards the first passage bridge of a stable L\'evy process with no
negative jumps and its height process.\\

In section 2, we introduce  conditioned Galton-Watson forests and
their related coding first passage bridge, height process and
contour process. Section 3 is devoted to recall the definitions of
real trees and L\'evy trees and to state a number of important
results related to these notions.

\section{Discrete trees and forests.}\label{disc}
In all the sequel, an element $u$ of $(\mathbb{N}^*)^n$ is written
as $u=(u_{1}, \ldots u_{n})$ and we set $|u|=n$. Let
\[\mathbb{U}=\bigcup_{n=0}^{\infty} (\mathbb{N}^*)^n,\]
where $\mathbb{N}^*=\{1, 2, \ldots\}$ and by convention
$(\mathbb{N}^*)^0=\{\emptyset\}$.  The concatenation of two elements
of $\mathbb{U}$, let us say $u=(u_{1}, \ldots u_{n})$ and $v=(v_{1},
\ldots,v_{m})$ is denoted by $uv=(u_{1}, \ldots u_{n},v_{1},
\ldots,v_{m})$.
 A discrete rooted tree is an element $\tau$ of the set
$\mathbb{U}$ which satisfies:
\begin{itemize}
\item[(i)] $\emptyset\in\tau$,
\item[(ii)] If $v\in\tau$ and $v=uj$ for some $j\in\mathbb{N}^*$,
then $u\in\tau$.
\item[(iii)] For every $u\in\tau$, there exists a number
$k_u(\tau)\ge0$, such that $uj\in\tau$ if and only if $1\le j\le
k_u(\tau)$.
\end{itemize}
In this definition, $k_u(\tau)$ represents the number of children of
the vertex $u$. We denote by $\mathbb{T}$ the set of all rooted
trees. The total cardinality of an element $\tau\in\mathbb{T}$ will
be denoted by $\zeta(\tau)$, (we emphasize that the root is counted
in $\zeta(\tau)$). If $\tau\in\mathbb{T}$ and $u\in\tau$, then we
define the shifted tree at the vertex $u$ by
\[\theta_u(\tau)=\{v\in\mathbb{U}:uv\in\tau\}\,.\]
We say that $u\in \tau$ is a leaf of $\tau$ if $k_{u}(\tau)=0$.

 Then we consider a probability measure $\mu$ on $\mathbb{Z}_+$, such that
\[\sum_{k=0}^\infty k\mu(k)\le
1\;\;\;\;\mbox{and}\;\;\;\;\mu(1)<1\,.
\]
\noindent The law of the Galton-Watson tree with offspring
distribution $\mu$ is the unique probability measure
$\mathbb{Q}_{\mu}$ on $\mathbb{T}$ such that:
\begin{itemize}
\item[$(i)$] $\mathbb{Q}_\mu(k_\emptyset=j)=\mu(j)$,
$j\in\mathbb{Z}_+$.
\item[$(ii)$] For every $j\ge1$, with $\mu(j)>0$, the shifted trees
$\theta_1(\tau),\dots,\theta_j(\tau)$ are independent under the
conditional distribution $\mathbb{Q}_\mu(\,\cdot\,|\,k_\emptyset=j)$
and their conditional law is $\mathbb{Q}_\mu$.
\end{itemize}
A Galton-Watson forest with offspring distribution $\mu$ is a finite
or infinite sequence of independent Galton-Watson trees with
offspring distribution $\mu$. It will be denoted by ${\cal
F}=(\tau_k)$. With a misuse of notation, we will denote by
$\mathbb{Q}_\mu$ the law on $(\mathbb{T})^{\mathbb{N}^*}$ of a
Galton-Watson forest with offspring distribution $\mu$.\\
It is known that the G-W process associated to a G-W tree or forest
does not code entirely its genealogy. In the aim of doing so, other
(coding) real valued processes have been defined. Amongst such
processes one can cite the {\it contour process}, the {\it height
process} and the associated random walk which will be called here
{\it the coding walk} and which is sometimes referred to as the
Luckazievicks path.
\begin{defi}
We denote by
$u_\tau(0)=\emptyset,\,u_\tau(1)=1,\dots,\,u_\tau(\zeta-1)$ the
elements of a tree $\tau$ which are enumerated in the
lexicographical order $($when no confusion is possible, we will
simply write $u(n)$ for $u_\tau(n)$$)$. Let us denote by $|u(n)|$
the rank of the generation of a vertex $u(n)\in\tau$.
\begin{itemize}
\item[$(1)$] The height function of a tree $\tau$ is defined by
\[n\mapsto H_n(\tau)=|u(n)|,\,0\le n\le \zeta(\tau)-1\,.\]
\item[$(2)$] The height function of a forest ${\cal
F}=(\tau_k)$ is defined by
\begin{eqnarray*} &&n\mapsto H_n({\cal
F})=H_{n-(\zeta(\tau_0)+\dots+\zeta(\tau_{k-1}))}(\tau_k),\\
&&\;\;\;\;\;\;\;\;\;\;\;\;\;\;\;\;\mbox{if
$\;\zeta(\tau_0)+\dots+\zeta(\tau_{k-1})\le
n\le\zeta(\tau_0)+\dots+\zeta(\tau_k)-1$,}
\end{eqnarray*}
for $k\ge1$, and with the convention that $\zeta(\tau_0)=0$. If
there is a finite number of trees in the forest, say $j$,  then we
set $H_n({\cal F})=0$, for $n\ge \zeta(\tau_0)+\dots+\zeta(\tau_j)$.
\end{itemize}
\end{defi}
\noindent For two vertices $u$ and $v$  of a tree $\tau$, the
distance $d_\tau(u,v)$ is the number of edges of the unique
elementary path from $u$ to $v$. The height function may be
presented in a natural way as the distance between the visited
vertex and the root $\emptyset$, i.e.
$H_n(\tau)=d_\tau(\emptyset,u(n))$. Then we may check the following
relation
\begin{equation}\label{dist}
d_\tau(u(n),u(m))=H_n(\tau)+H_m(\tau)-2H_{k(n,m)}(\tau)\,,
\end{equation}
where $k(n,m)$ is the index of the last common ancestor of $u(n)$
and $u(m)$. It is not difficult to see that the height process of a
tree (resp. a forest) allows us to recover the entire structure of
this tree (resp. this forest). We say that it {\it codes} the
genealogy of the tree or the forest. Although this process is
natural and simple to define, its law is rather complicated to
characterize. In particular, $H$ is neither a Markov process nor a
martingale.

The contour process gives another characterization of the tree which
is easier to visualize. We suppose that the tree is embedded in a
half-plane in such a way that edges have length one. Informally, we
imagine the motion of a particle that starts at time $0$ from the
root of the tree and then explores the tree from the left to the
right continuously along each edge of $\tau$ at unit speed until all
edges have been explored and the particle has come back to the root.
Note that if $u(n)$ is a leaf, then the particle goes to $u(n+1)$,
taking the shortest way that consists first to move backward on the
line of descent from $u(n)$ to their last common ancestor
$u(n)\wedge u(n+1)$ and then to move forward along the single edge
between $u(n)\wedge u(n+1)$ to $u(n+1)$. Since it is clear that each
edge will be crossed twice, the total time needed to explore the
tree is $2(\zeta(\tau)-1)$. The value $C_t(\tau)$ of the contour
function at time $t\in [0,2(\zeta(\tau)-1)]$ is the distance (on the
continuous tree) between the position of the particle at time $s$
and the root. More precisely, let us denote by $l_1<l_2<\cdots<l_p$
the $p$ leaves of $\tau$ listed in lexicographical order. The
contour function $(C_t(\tau), 0\leq t\leq 2(\zeta(\tau)-1))$ is the
piecewise linear continuous path with slope equal to +1 or -1, that
takes successive local extremes with values: $0, |l_1|,|l_1\land
l_2|, |l_2|, \ldots |l_{p-1}\land l_p|, |l_p|$ and $0$. Then we set
$C_t(\tau)=0$, for $t\in[2(\zeta(\tau)-1),2\zeta(\tau)]$. It is
clear that $C(\tau)$ codes the genealogy of  $\tau$.

The contour process for a forest $\mathcal{F}=(\tau_{k})$ is the
concatenation of the processes $C(\tau_1),\dots,C(\tau_k),\dots$,
i.e. for $k\ge1$:
\[
C_{t}(\mathcal{F})=C_{t-2(\zeta(\tau_0)+\cdots+\zeta(\tau_{k-1}))}(\tau_k),
\;\textrm{ if }\; 2(\zeta(\tau_0)+\cdots+\zeta(\tau_{k-1}))\le t\le
2(\zeta(\tau_0)+\cdots+\zeta(\tau_{k})).
\]
If there is a finite number of trees, say $j$, in the forest, we set
$C_t(\mathcal{F})=0$, for
$t\ge2(\zeta(\tau_0)+\cdots+\zeta(\tau_{j}))$. Note that for each
tree $\tau_k$, $[2(\zeta(\tau_k)-1),2\zeta(\tau_k)]$ is the only
non-trivial subinterval of $[0,2\zeta(\tau_k)]$ on which $C(\tau_k)$
vanishes. This convention ensures that the contour process
$C(\mathcal{F})$ also codes the genealogy of the forest. However, it
has no "good properties" in law either.

In order to define a coding process whose law can easily be
described, most of the authors introduce the coding random walk
$S(\tau)$ which is defined as follows:
\[S_0=0\,,\;\;\;\;S_{n+1}(\tau)-S_n(\tau)=k_{u(n)}(\tau)-1,\;\;\;0\le n\le
\zeta(\tau)-1\,.\] Here again it is not very difficult to see that
the process $S(\tau)$ codes the genealogy of the tree $\tau$.
However, its construction requires a little bit more care than this
of $H(\tau)$ or $C(\tau)$. For each $n$, $S_n(\tau)$ is the sum of
all the younger brother of each of the ancestor of $u(n)$ including
$u(n)$ itself. For a forest ${\cal F}=(\tau_k)$, the process
$S({\cal F})$ is the concatenation of
$S(\tau_1),\dots,S(\tau_k),\dots$:
\begin{eqnarray*}
&&S_n({\cal F})=S_{n-(\zeta(\tau_0)+\dots+\zeta(\tau_{k-1}))}(\tau_k)-k+1,\\
&&\;\;\;\;\;\;\;\;\;\;\;\;\;\;\;\;\mbox{if
$\;\zeta(\tau_0)+\dots+\zeta(\tau_{k-1})\le
n\le\zeta(\tau_0)+\dots+\zeta(\tau_k)$.}
\end{eqnarray*}
If there is a finite number of trees $j$, then we set $S_n({\cal
F})=S_{\zeta(\tau_0)+\dots+\zeta(\tau_j)}({\cal F})$, for
$n\ge\zeta(\tau_0)+\dots+\zeta(\tau_j)$. From the construction of
$S(\tau_1)$ it appears that $S(\tau_1)$ is a random walk with
initial value $S_{0}=0$ and step distribution $\nu(k)=\mu(k+1)$,
$k=-1, 0, 1, \ldots$ which is killed when it first enters into the
negative half-line. Hence, when the number of trees is infinite,
 $S(\mathcal{F})$ is a downward skip free random walk on $\mathbb{Z}$
 with the law described above.\\

\begin{center}
\unitlength=1.1pt
\begin{picture}(350,140)

\thicklines \put(50,0){\line(0,1){30}} \thicklines
\put(50,0){\line(-1,1){30}} \thicklines \put(50,0){\line(1,1){30}}

\thicklines \put(20,30){\line(-1,3){15}} \thicklines
\put(20,30){\line(1,3){15}}

\thicklines \put(5,75){\line(0,1){30}} \thicklines
\put(5,75){\line(-1,3){10}} \thicklines \put(5,75){\line(1,3){10}}

\thicklines \put(35,75){\line(-1,3){10}} \thicklines
\put(35,75){\line(1,3){10}}

\put(46,-10){$\varnothing$} \put(13,30){\scriptsize 1}
\put(43,30){\scriptsize 9}\put(70,30){\scriptsize
10}\put(-6,75){\scriptsize 2}\put(27,75){\scriptsize 6}
\put(60,78){\scriptsize 11}\put(92,78){\scriptsize 12}

\put(-8,109){\scriptsize 3}\put(4,109){\scriptsize 4}
\put(14,109){\scriptsize 5}\put(24,109){\scriptsize 7}
\put(44,109){\scriptsize 8}

 \thicklines \put(80,30){\line(-1,3){15}}
\thicklines \put(80,30){\line(1,3){15}}

\put(25,-45){Rooted tree $\tau$}

\hspace*{1.3in}

\thinlines \put(100,0){\line(1,0){180}} \thinlines
\put(100,-30){\line(0,1){170}}

\put(100,0){\circle*{3}} \put(113,50){\circle*{3}}

\put(126,75){\circle*{3}}\put(139,125){\circle*{3}}\put(152,100){\circle*{3}}
\put(165,75){\circle*{3}}\put(178,50){\circle*{3}}\put(191,75){\circle*{3}}
\put(204,50){\circle*{3}}\put(217,25){\circle*{3}}\put(230,0){\circle*{3}}
\put(243,25){\circle*{3}}\put(256,0){\circle*{3}}\put(269,-25){\circle*{3}}

\thinlines \put(113,0){\line(0,1){2}} \thinlines
\put(126,0){\line(0,1){2}} \thinlines \put(139,0){\line(0,1){2}}
\put(152,0){\line(0,1){2}}\put(165,0){\line(0,1){2}}\put(178,0){\line(0,1){2}}
\put(191,0){\line(0,1){2}}\put(204,0){\line(0,1){2}}\put(217,0){\line(0,1){2}}
\put(230,0){\line(0,1){2}}\put(243,0){\line(0,1){2}}\put(256,0){\line(0,1){2}}
\put(269,0){\line(0,1){2}}

\thinlines \put(100,-25){\line(1,0){2}}\thinlines
\put(100,25){\line(1,0){2}} \thinlines \put(100,50){\line(1,0){2}}
\thinlines \put(100,75){\line(1,0){2}} \thinlines
\put(100,100){\line(1,0){2}} \thinlines \put(100,125){\line(1,0){2}}

\put(150,-45){Coding walk $S(\tau)$}

\put(112,-7){\scriptsize 1} \put(125,-7){\scriptsize 2}
\put(138,-7){\scriptsize 3} \put(151,-7){\scriptsize 4}
\put(164,-7){\scriptsize 5}\put(177,-7){\scriptsize 6}
\put(190,-7){\scriptsize 7}\put(203,-7){\scriptsize 8}
\put(216,-7){\scriptsize 9}\put(226,-7){\scriptsize 10}
\put(240,-7){\scriptsize 11}\put(254,-7){\scriptsize 12}

\put(93,-26){\scriptsize -1}\put(95,24){\scriptsize 1}
\put(95,49){\scriptsize 2} \put(95,74){\scriptsize 3}
\put(95,99){\scriptsize 4} \put(95,124){\scriptsize 5}

\end{picture}

\vskip 25mm

Figure 1

\end{center}

\vspace*{0.3cm}

\noindent Let us denote $H({\cal F})$, $C({\cal F})$ and $S({\cal
F})$ respectively by $H$, $C$ and $S$ when no confusion is possible.
In the sequel, we will have to use some path relationships between
$H$, $C$ and $S$ which we recall now. Let us suppose that ${\cal F}$
is infinite. It is established for instance in  \cite{DL,LL} that
\begin{equation}\label{sh}
H_n=\mbox{\rm card}\,\big\{0\le k\le n-1:S_k=\inf_{k\le j\le
n}S_j\big\}\,.
\end{equation}
This identity means that the height process at each time $n$ can be
interpreted as the amount of time that the random walk $S$ spends at
its future minimum before $n$. The following relationship between
$H$ and $C$ is stated in \cite{DL}: set $K_n=2n-H_n$, then
\begin{equation}\label{ch}
C_t=\left\{\begin{array}{ll} (H_n-(t-K_n))^+\,,&\mbox{if
$t\in[K_n,K_{n+1}-1]$}\\
(H_{n+1}-(K_{n+1}-t))^+\,,&\mbox{if $[K_{n+1}-1,K_{n+1}]$}\,.
\end{array}\right.
\end{equation}

\noindent For any integer $k\ge1$, we  denote by ${\cal F}^{k,n}$ a
G-W forest with $k$ trees conditioned to have $n$ vertices, that is
a forest with the same law as ${\cal F}=(\tau_1,\dots,\tau_k)$ under
the conditional law
$\mathbb{Q}_\mu(\,\cdot\,|\,\zeta(\tau_1)+\cdots+\zeta(\tau_k)=n)$.
The starting point of our work is the observation that ${\cal
F}^{k,n}$ can be coded by a downward skip free random walk
conditioned to first reach $-k$ at time $n$. An interpretation of
this result may be found in \cite{Pi}, Lemma 6.3 for instance.

\begin{proposition}\label{p2}
Let ${\cal F}=(\tau_j)$ be an infinite forest with offspring
distribution $\mu$ and $S$, $H$ and $C$ be respectively its coding
walk, its height process and its contour process. Let $W$ be a
random walk defined on a probability space $(\Omega,{\cal F},P)$
with the same law as $S$. We define $T_i^W=\inf\{j:W_j=-i\}$, for
$i\ge1$. Take $k$ and $n$ such that $P(T_k^W=n)>0$. Then under the
conditional law
$\mathbb{Q}_\mu(\,\cdot\,|\,\zeta(\tau_1)+\cdots+\zeta(\tau_k)=n)$,
\begin{itemize}
\item[$(1)$] The process $(S_j,\,0\le
j\le\zeta(\tau_1)+\dots+\zeta(\tau_k))$ has the same law as
$(W_j,\,0\le j\le T_k^W)$.
\end{itemize}
Moreover, define the processes $H_n^W=\mbox{\rm
card}\,\big\{k\in\{0,\dots,n-1\}:W_k=\inf_{k\le j\le n}W_j\big\}$
and $C^W$ using the height process $H^W$ as in $(\ref{ch})$, then
\begin{itemize}
 \item[$(2)$] the process
$(H_j,\,0\le j\le\zeta(\tau_1)+\dots+\zeta(\tau_k))$ has the same
law as the process $(H_j^W,\,0\le j\le T_k^W)$.
\end{itemize}

\begin{itemize}
 \item[$(3)$] the process
$(C_t,\,0\le 0\leq t \le 2(\zeta(\tau_1)+\dots+\zeta(\tau_k)))$ has
the same law as the process $(C_t^W,\,0\le t\le 2T_k^W)$.
\end{itemize}
\end{proposition}

\noindent It is also straightforward  that the  identities in law
involving separately the processes $H$, $S$ and $C$ in the above
proposition also hold for the triple $(H,S,C)$.  In the figure
below, we have represented an occurrence of the forest ${\cal
F}^{k,n}$ and its associated {\it coding first passage bridge}.
\begin{center}
\unitlength=1.1pt
\begin{picture}(350,140)

\thicklines \put(30,75){\line(0,1){12}} \thicklines
\put(30,75){\line(-1,1){12}} \thicklines \put(30,75){\line(1,1){12}}
\thicklines \put(18,87){\line(-1,3){4}} \thicklines
\put(18,87){\line(1,3){4}}\thicklines \put(42,87){\line(0,1){12}}
\thicklines \put(42,97){\line(0,1){10}} \thicklines
\put(42,97){\line(-1,2){5}} \thicklines \put(42,97){\line(1,2){5}}

\thicklines \put(65,75){\line(-1,2){6}} \thicklines
\put(65,75){\line(1,2){6}} \thicklines \put(71,87){\line(-1,2){5}}
\thicklines \put(71,87){\line(1,2){5}}

\thicklines \put(100,75){\line(0,1){10}} \thicklines
\put(100,85){\line(0,1){10}} \thicklines
\put(100,85){\line(-1,2){5}} \thicklines \put(100,85){\line(1,2){5}}
\thicklines \put(95,95){\line(-1,3){4}} \thicklines
\put(95,95){\line(1,3){4}}

\put(135,75){\circle*{2}}\put(170,75){\circle*{2}}
\put(205,75){\circle*{2}}\put(240,75){\circle*{2}}

\thicklines \put(310,75){\line(-1,2){6}} \thicklines
\put(310,75){\line(1,2){6}} \thicklines \put(316,87){\line(-1,2){5}}
\thicklines \put(316,87){\line(1,2){5}} \thicklines
\put(311,97){\line(0,1){10}} \thicklines
\put(311,97){\line(-1,2){5}} \thicklines \put(311,97){\line(1,2){5}}

\thicklines \put(275,75){\line(0,1){10}} \thicklines
\put(275,85){\line(-1,2){5}} \thicklines \put(275,85){\line(1,2){5}}
\thicklines \put(280,95){\line(-1,3){4}} \thicklines
\put(280,95){\line(1,3){4}} \thicklines
\put(276,107){\line(0,1){10}}

\put(28,60){$\tau_1$}\put(63,60){$\tau_2$}
\put(98,60){$\tau_3$}\put(273,60){$\tau_{k-1}$}
\put(308,60){$\tau_k$}

\put(29,70){\tiny 0}\put(13,85){\tiny 1}\put(11,102){\tiny
2}\put(21,102){\tiny 3}\put(26,85){\tiny 4}\put(36,85){\tiny 5}

\put(314,110){\tiny n-1}\put(320,100){\tiny $n$}

 \vspace*{0.3in}
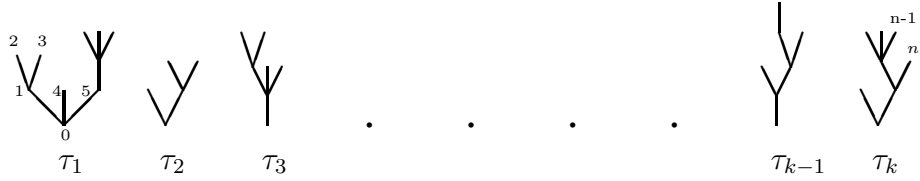
\put(55,40){\small Conditioned random forest: $k$ trees, $n$
vertices.}

\hspace*{0.4in}

\thinlines \put(0,-38){\line(1,0){300}} \thinlines
\put(10,-115){\line(0,1){145}}

\put(10,-38){\circle*{2}} \put(20,-18){\circle*{2}}
\put(30,-8){\circle*{2}}
\put(40,-18){\circle*{2}}\put(50,-28){\circle*{2}}\put(60,-38){\circle*{2}}
\put(70,-38){\circle*{2}}\put(80,-18){\circle*{2}}\put(90,-28){\circle*{2}}
\put(180,-68){\circle*{2}}\put(190,-68){\circle*{2}}\put(200,-78){\circle*{2}}

\put(210,-88){\circle*{2}}\put(220,-78){\circle*{2}}
\put(230,-88){\circle*{2}}\put(240,-98){\circle*{2}}\put(250,-78){\circle*{2}}
\put(260,-88){\circle*{2}}\put(270,-88){\circle*{2}}\put(280,-98){\circle*{2}}
\put(290,-108){\circle*{2}}

\thinlines \put(20,-40){\line(0,1){2}}
\put(30,-40){\line(0,1){2}}\put(40,-40){\line(0,1){2}}\put(50,-40){\line(0,1){2}}
\put(60,-40){\line(0,1){2}}\put(70,-40){\line(0,1){2}}\put(80,-40){\line(0,1){2}}
\put(90,-40){\line(0,1){2}}

\put(180,-40){\line(0,1){2}}\put(190,-40){\line(0,1){2}}\put(200,-40){\line(0,1){2}}
\put(210,-40){\line(0,1){2}}\put(220,-40){\line(0,1){2}}
\put(230,-40){\line(0,1){2}}\put(240,-40){\line(0,1){2}}\put(250,-40){\line(0,1){2}}
\put(260,-40){\line(0,1){2}}\put(270,-40){\line(0,1){2}}\put(280,-40){\line(0,1){2}}
\put(290,-40){\line(0,1){2}}

\thinlines \put(8,22){\line(1,0){2}} \thinlines
\put(8,12){\line(1,0){2}}\thinlines \put(8,2){\line(1,0){2}}
\thinlines \put(8,-8){\line(1,0){2}}\thinlines
\put(8,-18){\line(1,0){2}}\thinlines
\put(8,-28){\line(1,0){2}}\thinlines \put(8,-38){\line(1,0){2}}
\thinlines \put(8,-48){\line(1,0){2}}\thinlines
\put(8,-58){\line(1,0){2}} \thinlines
\put(8,-68){\line(1,0){2}}\thinlines \put(8,-78){\line(1,0){2}}
\thinlines\put(8,-88){\line(1,0){2}} \thinlines
\put(8,-98){\line(1,0){2}} \thinlines \put(8,-108){\line(1,0){2}}

\put(18,-47){{\tiny $1$}} \put(28,-47){{\tiny $2$}}
\put(285,-47){{\tiny $n$}}

\put(3,-36){{\tiny $0$}} \put(-3,-50){{\tiny $-1$}}
\put(-3,-60){{\tiny $-2$}} \put(-3,-110){{\tiny $-k$}}

\thinlines \put(10,-108){\line(1,0){8}} \put(25,-108){\line(1,0){8}}
\put(40,-108){\line(1,0){8}} \put(55,-108){\line(1,0){8}}
\put(70,-108){\line(1,0){8}} \put(85,-108){\line(1,0){8}}
\put(100,-108){\line(1,0){8}} \put(115,-108){\line(1,0){8}}
\put(130,-108){\line(1,0){8}} \put(145,-108){\line(1,0){8}}
\put(160,-108){\line(1,0){8}} \put(175,-108){\line(1,0){8}}
\put(190,-108){\line(1,0){8}} \put(205,-108){\line(1,0){8}}
\put(220,-108){\line(1,0){8}} \put(235,-108){\line(1,0){8}}
\put(250,-108){\line(1,0){8}} \put(265,-108){\line(1,0){8}}
\put(281,-108){\line(1,0){9}}

\thinlines \put(290,-108){\line(0,1){8}}
\put(290,-95){\line(0,1){8}} \put(290,-80){\line(0,1){8}}
\put(290,-65){\line(0,1){8}} \put(290,-48){\line(0,1){8}}

\put(80,-130){\small Coding first passage bridge}

\end{picture}

\vskip 55mm

Figure 2

\end{center}
 \vspace*{0.3in}
In Section \ref{constforest}, we will present a continuous time
version of this result, but before we need to introduce the
continuous time setting of L\'evy trees and forests.

\section{Coding real trees and forests}\label{codcont}

Discrete trees may be considered in an obvious way as compact metric
spaces with no loops. Such metric spaces are special cases of
$\R$-trees which are defined hereafter. Similarly to the discrete
case, an $\R$-forest is any collection of $\R$-trees. In this
section we keep the same notations as in Duquesne and Le Gall's
articles \cite{DL} and \cite{DL2}. The following formal definition
of $\R$-trees is now classical and originates from $T$-theory. It
may be found for instance in \cite{DMT}.
\begin{defi}
A metric space $({\cal T},d)$ is an $\R$-tree if for every
$\sigma_1,\sigma_2\in{\cal T}$,
\begin{itemize}
\item[$1.$] There is a unique map $f_{\sigma_1,\sigma_2}$ from
$[0,d(\sigma_1,\sigma_2)]$ into ${\cal T}$ such that
$f_{\sigma_1,\sigma_2}(0)=\sigma_1$ and
$f_{\sigma_1,\sigma_2}(d(\sigma_1,\sigma_2))=\sigma_2$.
\item[$2.$] If  $g$ is a continuous injective map from $[0,1]$ into
${\cal T}$ such that $g(0)=\sigma_1$ and $g(1)=\sigma_2$, we have
\[g([0,1])=f_{\sigma_1,\sigma_2}([0,d(\sigma_1,\sigma_2)])\,.\]
\end{itemize}
A rooted $\R$-tree is an $\R$-tree $({\cal T},d)$ with a
distinguished vertex $\rho=\rho({\cal T})$ called the root. An
$\R$-forest is any collection of rooted $\R$-trees: ${\cal
F}=\{({\cal T}_i,d_i),\,i\in{\cal I}\}$.
\end{defi}
\noindent A construction of some particular cases of such metric
spaces has been given by Aldous \cite{Al} and is described in
\cite{DL2} in a more general setting. Let
$f:[0,\infty)\rightarrow[0,\infty)$ be a continuous function with
compact support, such that $f(0)=0$. For $0\le s\le t$, we define
\begin{equation}\label{distance}
d_f(s,t)=f(s)+f(t)-2\inf_{u\in[s,t]}f(u)
\end{equation}
and the equivalence relation by
\[s\sim t\;\;\;\mbox{if and only if}\;\;\;d_f(s,t)=0\,.\]
(Note that $d_f(s,t)=0$ if and only if
$f(s)=f(t)=\inf_{u\in[s,t]}f(u)$.) Then the projection of $d_f$ on
the quotient space
\[{\cal T}_f=[0,\infty)/\sim\]
defines a distance. This distance will also be denoted by $d_f$.
\begin{theorem}\label{tree}
The metric space $({\cal T}_f,d_f)$ is a compact $\R$-tree.
\end{theorem}
\noindent Denote by $p_f:[0,\infty)\rightarrow{\cal T}_f$ the
canonical projection. The vertex $\rho=p_f(0)$ will be chosen as the
root of ${\cal T}_f$. It has recently been proved by Duquesne
\cite{Du} that any $\R$-tree (satisfying some rather weak
assumptions) may by represented as $({\cal T}_f,d_f)$ where $f$ is a
left continuous function with right limits and without positive
jumps.

When no confusion is possible with the discrete case, the space of
$\R$-trees will also be denoted by $\mathbb{T}$. It is endowed with
the Gromov-Hausdorff distance, $d_{GH}$ which we briefly recall now.
For a metric space $(E,\delta)$ and $K$, $K'$ two subspaces of $E$,
$\delta_{\mbox{\tiny Haus}}(K,K')$ will denote the Hausdorff
distance between $K$ and $K'$. Then we define the distance between
${\cal T}$ and ${\cal T}'$ by:
\[d_{GH}({\cal T},{\cal
T}')=\inf\left(\delta_{\mbox{\tiny Haus}}(\varphi({\cal
T}),\varphi'({\cal
T}'))\vee\delta(\varphi(\rho),\varphi'(\rho'))\right)\,,\] where the
infimum is taken over all isometric embeddings $\varphi:{\cal
T}\rightarrow E$ and $\varphi':{\cal T}'\rightarrow E$ of ${\cal T}$
and ${\cal T}'$ into a common metric space $(E,\delta)$. We refer to
Chapter 3 of Evans \cite{Ev} and the references therein for a
complete description of the Gromov-Hausdorff topology. It is
important to note that the space $(\mathbb{T},d_{GH})$ is complete
and separable, see for instance Theorem 3.23 of \cite{Ev} or
\cite{EPW}.

In the remainder of this section, we will recall from \cite{DL2} the
definition of L\'evy trees and we state this of L\'evy forests. Let
$(\p_x)$, $x\in\R$ be a sequence of probability measures on the
Skorokhod space $\mathbb{D}$ of c\`adl\`ag paths from $[0,\infty)$
to $\R$ such that for each $x\in\R$, the canonical process  $X$ is a
L\'evy process with no negative jumps. Set $\p=\p_0$, so $\p_x$ is
the law of $X+x$ under $\p$. We suppose that the characteristic
exponent $\psi$ of $X$ (i.e. $\e(e^{-\lambda
X_t})=e^{t\psi(\lambda)}$, $\lambda\in\R$) satisfies the following
condition:
\begin{equation}\label{cond}
\int_1^{\infty}\frac{du}{\psi(u)}<\infty\,.\end{equation}  By
analogy with the discrete case, the continuous time height process
$\bar{H}$ is the measure (in a sense which is to be defined) of the
set $\{s\le t:X_s=\inf_{s\le r\le t}X_r\}$. A rigorous meaning to
this measure is given by the following result due to Le Jan and Le
Gall \cite{LL}, see also \cite{DL}. Define $I_t^s=\inf_{s\le u\le
t}X_u$. There is a sequence of positive real numbers
$(\varepsilon_k)$ which decreases to $0$ such that for any $t$, the
limit
\begin{equation}\label{5289}
\bar{H}_t\eqdef\lim_{k\rightarrow+\infty}\frac1{\varepsilon_k}\int_0^t
\ind_{\{X_s-I^s_t<\varepsilon_k\}}\,ds\end{equation} exists a.s. It
is also proved in \cite{LL} that under assumption (\ref{cond}),
$\bar{H}$ is a continuous process, so that each of its positive
excursion codes a real tree in the sense of Aldous. We easily deduce
from this definition that the height process $\bar{H}$ is a
functional of the L\'evy process reflected at its minimum, i.e.
$X-I$, where $I:=I^0$. In particular, when $\alpha=2$, $\bar{H}$ is
equal in law to the reflected process multiplied by a constant. It
is well known that $X-I$ is a strong Markov process. Moreover, under
our assumptions, 0 is regular for itself for this process and we can
check that the process $-I$ is a local time at level 0. We denote by
$N$ the corresponding It\^o measure of the excursions away from 0.

In order to define the L\'evy forest, we need to introduce the local
times of the height process $\bar{H}$. It is proved in \cite{DL}
that for any level $a\ge0$, there exists a continuous increasing
process $(L_t^a,\,t\ge0)$ which is defined by the approximation:
\begin{equation}\label{loctime}\lim_{\varepsilon\downarrow0}
\e\left(\sup_{0\le s\le
t}\left|\frac1\varepsilon\int_0^sdu\ind_{\{a<\bar{H}_u\le
a+\varepsilon\}}-L_s^a\right|\right)=0\,.\end{equation} The support
of the measure $dL_t^a$ is contained in the set
$\{t\ge0:\bar{H}_t=a\}$ and we readily check that $L^0=-I$. Then we
may define the Poisson  point process of the excursions away from 0
of the process $\bar{H}$ as follows. Let $T_u=\inf\{t:-I_t\ge u\}$
be the right continuous inverse of the local time at 0 of the
reflected process $X-I$ (or equivalently of $H$).  The time $T_u$
corresponds to the first passage time of $X$ bellow $-u$. Set
$T_{0-}=0$ and for all $u\ge0$,
\[e_u(v)=\left\{\begin{array}{ll}
\bar{H}_{T_{u-}+v}\,,&\mbox{if $0\le v\le T_u-T_{u-}$}\\
0\,,&\mbox{if $v> T_u-T_{u-}$}\end{array}\right.\,.\] For each
$u\ge0$, we may define the tree $({\cal T}_{e_u},d_{e_u})$ under
$\p$ as in the beginning of this section. We easily deduce from the
Markov property of $X-I$ that under the probability measure $\p$,
the process $\{({\cal T}_{e_u},d_{e_u}),u\ge0\}$ is a Poisson point
process whose characteristic measure is the law of the random real
tree $({\cal T}_{\bar{H}},d_{\bar{H}})$ under $N$. By analogy to the
discrete case, this Poisson point process, as a $\mathbb{T}$-valued
process, provides a natural definition for the L\'evy forest.

\begin{defi}
The L\'evy tree is the real tree $({\cal T}_{\bar{H}},d_{\bar{H}})$
coded by the function $\bar{H}$ under the measure $N$. We denote by
$\Theta(d{\cal T})$ the $\sigma$-finite measure on $\mathbb{T}$
which is the law of the L\'evy tree ${\cal T}_{\bar{H}}$ under $N$.
The L\'evy forest ${\cal F}_{\bar{H}}$ is the Poisson point process
\[({\cal F}_{\bar{H}}(u),\,u\ge0)\eqdef\{({\cal
T}_{e_u},d_{e_u}),u\ge0\}\] which has for characteristic measure
$\Theta(d{\cal T})$ under $\p$. For each $s>0$, the process ${\cal
F}_{\bar{H}}^s\eqdef\{({\cal T}_{e_u},d_{e_u}),0\le u\le s\}$ under
$\p$ will be called the L\'evy forest of {\it size} $s$.
\end{defi}
\noindent Such a definition of a L\'evy forest has already been
introduced in \cite{Pi}, Proposition 7.8 in the Brownian setting.
However in this work, it is observed that the Brownian forest may
also simply be defined as the real tree coded by the function
$\bar{H}$ under  law $\p$. We also refer to \cite{WP} where the
Brownian forest is understood in this way.
 Similarly, the L\'evy forest with size $s$ may  be defined as
the compact real tree coded by the continuous function with compact
support $(\bar{H}_u,\,0\le u\le T_s)$ under law $\p$. These
definitions are more natural when considering convergence of
sequences of real forests and we will make appeal to them in section
5, see Corollary \ref{cor}.

We will simply denote the L\'evy tree and the L\'evy forest
respectively  by ${\cal T}_{\bar{H}}$, ${\cal F}_{\bar{H}}$ or
${\cal F}_{\bar{H}}^s$, the corresponding distances  being implicit.
When $X$ is stable, condition (\ref{cond}) is satisfied if and only
if its index $\alpha$ satisfies $\alpha\in(1,2)$. We may check, as a
consequence of (\ref{5289}), that $\bar{H}$ is a self-similar
process with index $\alpha/(\alpha-1)$, i.e.:
\[(\bar{H}_t,\,t\ge0)\ed(k^{(\alpha-1)/\alpha}
\bar{H}_{kt},\,t\ge0)\,,\;\;\;\mbox{for all $k>0$.}\] In this case,
the L\'evy tree ${\cal T}_{\bar{H}}$ associated to the stable
mechanism is called the $\alpha$-stable L\'evy tree and its law is
denoted by $\Theta_\alpha(d{\cal T})$. This random metric space also
inherits from $X$ a scaling property which may be stated as follows:
for any $a>0$, we denote by $a{\cal T}_{\bar{H}}$ the L\'evy tree
${\cal T}_{\bar{H}}$ endowed with the distance $ad_{\bar{H}}$, i.e.
\begin{equation}\label{scalingtree}a{\cal T}_{\bar{H}}\eqdef
({\cal T}_{\bar{H}},a d_{\bar{H}})\,.\end{equation} Then the law of
$a{\cal T}_{\bar{H}}$ under $\Theta_\alpha(d{\cal T})$ is
$a^{\frac1{\alpha-1}}\Theta_\alpha(d{\cal T})$. This property is
stated in \cite{LG} Proposition 4.3 and \cite{DL3} where other
fractal properties of stable trees are considered.

\newpage

\section{Construction of the conditioned L\'evy
forest}\label{constforest}

In this section we present the continuous analogue of the forest
${\cal F}^{k,n}$ introduced in section 2. In particular, we define
the total mass of the L\'evy forest of a given size $s$. Then we
define the L\'evy forest of size $s$ conditioned by its total mass.
In the stable case, we give a construction of this conditioned
forest from the unconditioned
forest.\\

We begin with the definition of the measure $\ell^{a,u}$  which
represents a local time at level $a>0$ for the L\'evy tree ${\cal
T}_{e_u}$. For all  $a>0$, $u\ge0$ and for every bounded and
continuous function $\varphi$ on ${\cal T}_{e_u}$, the finite
measure $\ell^{a,u}$ is defined by:
\begin{equation}\label{local}\langle\ell^{a,u},\varphi\rangle=
\int_{0}^{T_u-T_{u-}}
dL_{T_u+v}^a\varphi(p_{e_u}(v))\,,\end{equation} where we recall
from the previous section that $p_{e_u}$ is the canonical projection
from $[0,\infty)$ onto ${\cal T}_{e_u}$ for the equivalence relation
$\sim$ and $(L^a_{u})$ is the local time at level $a$ of $\bar{H}$.
Then the mass measure of the L\'evy tree ${\cal T}_{e_u}$ is
\begin{equation}\label{mass}\mbox{\bf m}_{u}=\int_0^\infty
da\,l^{a,u}\end{equation} and the total mass of the tree is
 $\mbox{\bf m}_{u}({\cal T}_{e_u})$.
 Now we fix $s>0$; the total mass of the forest of size
$s$, ${\cal F}_{\bar{H}}^s$ is naturally given by
\[\mbox{\bf M}_s=\sum_{0\le u\le s}\mbox{\bf m}_{u}({\cal T}_{e_u})\,.\]
\begin{proposition}\label{prop2}
$\p$-almost surely $T_s=\mbox{\bf M}_s$.
\end{proposition}
\noindent {\it Proof}. It follows from the definitions (\ref{local})
and (\ref{mass}) that for each tree ${\cal T}_{e_u}$, the mass
measure $\mbox{\bf m}_{u}$ coincides with the image of the Lebesgue
measure on $[0,T_u-T_{u-}]$ under the mapping $v\mapsto p_{e_u}(v)$.
Thus, the total mass $\mbox{\bf m}_{u}({\cal T}_{e_u})$ of each tree
${\cal T}_{e_u}$ is $T_u-T_{u-}$. This implies the result.\QED
\noindent Then we will construct processes which encode the
genealogy of the L\'evy forest of size $s$ conditioned to have a
mass equal to $t>0$. From the analogy with the discrete case in
Proposition \ref{p2}, the natural candidates may informally be
defined as:
\begin{eqnarray*}
X^{br}&\eqdef&[(X_u,\,0\le u\le T_s)\,|\,T_s=t]\\
\bar{H}^{br}&\eqdef&[(\bar{H}_u,\,0\le u\le T_s)\,|\,T_s=t]\,.
\end{eqnarray*}
When $X$ is the Brownian motion, the process $X^{br}$ is called the
{\it first passage bridge}, see \cite{BCP}. In order to give a
proper definition in the general case, we need the additional
assumption:
\[\mbox{\it The semigroup of $(X,\p)$ is absolutely continuous with
respect to the Lebesgue measure.}\] Then denote by $p_t(\cdot)$ the
density of the semigroup of $X$, by ${\cal G}^X_u\eqdef
\sigma\{X_v,v\le u\}$, $u\ge0$ the $\sigma$-field generated by $X$
and set $\hat{p}_t(x)=p_t(-x)$.
\begin{lemma}
The probability measure which is defined on each ${\cal G}^X_u$ by
\begin{equation}\label{fpb}
\p(X^{br}\in \Lambda_u)=\e\left(\ind_{\{X\in\Lambda_u,\,u<T_{s}\}}
\frac{t(s+X_u)}{s(t-u)}\frac{\hat{p}_{t-u}(s+X_u)}{\hat{p}_t(s)}\right)
\,,\;\;u<t\,,\;\;\;\Lambda_u\in{\cal G}_u^X\,,
\end{equation}
is a regular version of the conditional law of $(X_u,\,0\le u\le
T_s)$ given $T_s=t$, in the sense that for all $u>0$, for
$\lambda$-a.e.~$s>0$ and $\lambda$-a.e.~$t>u$,
\[\p(X^{br}\in
\Lambda_u)=\lim_{\varepsilon\downarrow0}\p(X\in\Lambda_u\,|
\,|T_s-t|<\varepsilon)\,,\] where $\lambda$ is the Lebesgue measure.
\end{lemma}
\noindent {\it Proof}. Let $u<t$, $\Lambda_u\in{\cal G}_u^X$ and
$\varepsilon<t-u$. From the Markov property, we may write
\begin{eqnarray}\p(X\in\Lambda_u\,|
\,|T_s-t|<\varepsilon)&=&\e\left(\ind_{\{X\in\Lambda_u\}}\frac{\ind_{\{
\,|T_s-t|<\varepsilon\}}}{\p(|T_s-t|<\varepsilon)}\right)\nonumber\\
&=&\e\left(\ind_{\{X\in\Lambda_u,u<T_s\}}\frac{\p_{X_u}(
|T_s-(t-u)|<\varepsilon)}{\p(|T_s-t|<\varepsilon)}\right)\,.\label{23914}
\end{eqnarray}
On the other hand, from Corollary VII.3 in \cite{Be} one has,
\begin{equation}\label{espacetemps}
t\p(T_s\in dt)\,ds=s\hat{p}_t(s)\,dt\,ds\,.\end{equation} Hence, for
all $x\in\R$, for all $u>0$, for $\lambda$-a.e. $s>0$ and
$\lambda$-a.e. $t>u$,
\[\lim_{\varepsilon\downarrow0}\frac{\p_x(
|T_s-(t-u)|<\varepsilon)}{\p(|T_s-t|<\varepsilon)}=
\frac{t(s+x)}{s(t-u)}\frac{\hat{p}_{t-u}(s+x)}{\hat{p}_t(s)}\,.\]
Moreover we can check from (\ref{espacetemps}) that $\e\left(
\frac{t(s+X_u)}{s(t-u)}\frac{\hat{p}_{t-u}(s+X_u)}{\hat{p}_t(s)}\right)<+\infty$
for $\lambda$-a.e.~$t$, so the result follows from (\ref{23914}) and
Fatou's lemma.\QED

\noindent We may now construct a height process $\bar{H}^{br}$ from
the path of the first passage bridge $X^{br}$ exactly as $\bar{H}$
is constructed from $X$ in (\ref{5289}) or in Definition 1.2.1 of
\cite{DL} and check that the law of $\bar{H}^{br}$ is a regular
version of the conditional law of $(\bar{H}_u,0\le u\le T_s)$ given
$T_s=t$. Call $(e^{s,t}_u,\,0\le u\le s)$ the excursion process of
$\bar{H}^{br}$, that is in particular
\[\mbox{$(e^{s,t}_u,\,0\le u\le s)$ has the same law as $(e_u,\,0\le u\le
s)$ given $T_s=t$}\,.\] The following proposition is a
straightforward consequence of the above definition and Proposition
\ref{prop2}.
\begin{proposition}\label{prop4}
The law of the process $\{({\cal T}_{ e^{s,t}_u},d_{{
e^{s,t}_u}}),\,0\le v\le s\}$ is a regular version of the law of the
forest of size $s$, ${\cal F}_{\bar{H}}^s$ given $\mbox{\bf M}_s=t$.
\end{proposition}
\noindent We will denote by $({\cal F}^{s,t}_{\bar{H}}(u),\,0\le
u\le s)$ a process with values in $\mathbb{T}$ whose law under $\p$
is this of the L\'evy forest of size $s$ conditioned by $\mbox{\bf
M}_s=t$, i.e. conditioned
 to have a mass equal to $t$.\\

In the remainder of this section, we will consider the case when the
driving L\'evy process is stable. We suppose that its index $\alpha$
belongs to $(1,2]$ so that condition (\ref{cond}) is satisfied. We
will give a pathwise construction of the processes
$(X^{br},\bar{H}^{br})$ from the path of the original processes
$(X,\bar{H})$. This result leads to the following realization of the
L\'evy forest of size $s$ conditioned by its mass. From now on, with
no loss of generality, we suppose that $t=1$.

\begin{theorem}\label{constforet} Define $g=\sup\{u\le 1:
T_{u^{1/\alpha}}=s\cdot u\}$. \begin{itemize} \item[$(1)$]
$\p$-almost surely,
\[0<g<1\,.\]
\item[$(2)$] Under $\p$, the rescaled process
\begin{equation}\label{54927}
(g^{(1-\alpha)/\alpha}\bar{H}(gu),\,0\le u\le 1) \end{equation} has
the same law as $\bar{H}^{br}$ and is independent of $g$.
\item[$(3)$] The forest ${\cal F}^{s,1}_{\bar{H}}$ of size $s$ and mass
$1$ may be constructed from the rescaled process defined in
$(\ref{54927})$, i.e. if we denote by $u\mapsto \epsilon_u\eqdef
(g^{(1-\alpha)/\alpha}e_u(gv),v\ge 0)$ its process of excursions
away from $0$, then under $\p$,  ${\cal
F}_{\bar{H}}^{s,1}\ed\{({\cal T}_{ \epsilon_u},d_{{
\epsilon_u}}),\,0\le u\le s\}$.
\end{itemize}
\end{theorem}
\noindent {\it Proof}. The process $T_u=\inf\{v:I_v\le -u\}$ is a
stable subordinator with index $1/\alpha$. Therefore,
\[T_u<su^{\alpha}\,,\;\;\;\mbox{ i.o. as $u\downarrow0\;\;\;$
and}\;\;\;\;T_u>su^{\alpha}\,,\;\;\;\mbox{ i.o. as
$u\downarrow0$.}\] Indeed, if $u_n\downarrow0$ then
$\p(T_{u_n}<su_n^{\alpha})=\p(T_1<s)>0$, so that
$\p(\limsup_n\{T_{u_n}<su_n^{\alpha}\})\ge\p(T_1<s)>0$. But $T$
satisfies Blumenthal 0-1 law, so this probability is 1. The same
arguments prove that $\p(\limsup_n\{T_{u_n}>su_n^{\alpha}\})=1$ for
any sequence $u_n\downarrow0$. Since $T$ has only positive jumps, we
deduce that $T_u=su^{\alpha}$ infinitely often as $u$ tends to 0, so
we have proved the first part of the theorem.

The rest of the proof is a consequence of the following lemma.
\begin{lemma}\label{6581}
The first passage bridge $X^{br}$ enjoys the following path
construction:
\[X^{br}\ed (g^{-1/\alpha}X(gu),\,0\le u\le 1)\,.\]
Moreover, the process $(g^{-1/\alpha}X(gu),\,0\le u\le 1)$ is
independent of $g$.
\end{lemma}
\noindent {\it Proof}.  First note that for any $t>0$ the bivariate
random variable $(X_t,I_t)$ under $\p$ is absolutely continuous with
respect to the Lebesgue measure and there is a version of its
density which is continuous. Indeed from the Markov property and
(\ref{espacetemps}), one has for all $x\in\R$ and $y\ge0$,
\begin{eqnarray*}
\p(I_t\le y\,|\,X_t=x)&=&\e\left(\ind_{\{T_y\le
t\}}\frac{p_{t-T_y}(x-y)}{p_t(0)}\right)\\
&=&\int_0^t\frac ys\hat{p}_s(y)\frac{p_{t-s}(x-y)}{p_t(0)}\,ds\,.
\end{eqnarray*}
Looking at the expressions of $\hat{p}_t(x)$ and $p_t(x)$ obtained
from the Fourier inverse of the characteristic exponent of $X$ and
$-X$ respectively, we see that theses functions are continuously
differentiable and that their derivatives are continuous in $t$. It
allows us to conclude.

Now let us consider the two dimensional self-similar strong Markov
process $Y\eqdef(X,I)$ with state space $\{(x,y)\in\R^2:y\le x\}$.
>From our preceding remark, the semi-group
$q_t((x,y),(dx',dy'))=\p(X_t+x\in dx',y\wedge(I_t+x)\in dy')$ of $Y$
is absolutely continuous with respect to the Lebesgue measure and
there is a version of its density which is continuous. Denote by
$q_t((x,y),(x',y'))$  this version. We derive from
(\ref{espacetemps}) that for all $-s\le x$,
\begin{equation}\label{azv}
q_t((x,y),(-s,-s))=\ind_{\{y\ge
-s\}}\frac1t\hat{p}_t(s+x)\,.\end{equation}
 Then we may apply a
result due to Fitzsimmons, Pitman and Yor \cite{FPY} which asserts
that the inhomogenous Markov process on $[0,t]$, whose law is
defined by
\begin{equation}\label{bridge}
\e\left(H(Y_u,v\le
u)\frac{q_{t-u}(Y_u,(x',y'))}{q_t((x,y),(x',y'))}\,|\,Y_0=(x,y)
\right)\,,\;\;\;0\le u<t\,,\end{equation} where $H$ is a measurable
functional on $C([0,u],\R^2)$, is a regular version of the
conditional law of $(Y_v,\,0\le v\le t)$ given $Y_t=(x',y')$, under
$\p(\,\cdot\,|\,Y_0=(x,y))$. This law is called the law of the
bridge of $Y$ from $(x,y)$ to $(x',y')$ with length $t$. Then from
(\ref{azv}), the law which is defined in (\ref{bridge}), when
specifying it on the first coordinate and for $(x,y)=(0,0)$ and
$(x',y')=(-s,-s)$, corresponds to the law of the first passage
bridge which is defined in (\ref{fpb}).

It remains to apply another result which may also be found in
\cite{FPY}: observe that $g$ is a backward time for $Y$ in the sense
of \cite{FPY}. Indeed $g$ we may check that
$g=\sup\{u\le1:X_u=-su^{1/\alpha},\,X_u=I_u\}$, so that for all
$u>0$, $\{g>u\}\in\sigma(Y_v:v\ge u)$. Then from Corollary 3 in
\cite{FPY}, conditionally on $g$, the process $(Y_u,\,0\le u\le g)$
under $\p(\,\cdot\,|\,Y_0=(0,0))$ has the law of a bridge from
$(0,0)$ to $Y_g$ with length $g$. (This result has been obtained and
studied in a greater generality in \cite{CU}.) But from the
definition of $g$, we have $Y_g=(-sg^{1/\alpha},-sg^{1/\alpha})$, so
 from the self-similarity of $Y$, under $\p$ the process
\[({g^{-1/\alpha}}Y(g\cdot u)\,,0\le u\le 1)\]
has the law of the bridge of $Y$ from $(0,0)$ to $(-s,-s)$ with
length $1$. The lemma follows by specifying this result on the first
coordinate.\QED

\noindent The  second part of the theorem is a consequence of  Lemma
\ref{6581}, the construction of $\bar{H}^{br}$ from $X^{br}$  and
 the scaling property of $\bar{H}$. The third part follows from the definition
of the conditioned forest ${\cal F}^{s,1}_{\bar{H}}$ in Proposition
\ref{prop2} and the second part of this theorem.\QED

\newpage

\section{Invariance principles}\label{invprinc}

We know from Lamperti that the only possible limits of sequences of
re-scaled G-W processes are continuous state branching processes.
Then a question which arises is: when can we say that the whole
genealogy of the tree or the forest converges~? In particular, do
the height process, the contour process and the coding walk converge
after a suitable re-scaling~? This question has now been completely
solved by Duquesne and Le Gall \cite{DL}. Then one may ask the same
for the trees or forests conditioned by their mass. In \cite{Du},
Duquesne proved that when the law $\nu$ is in the domain of
attraction of a stable law, the height process, the contour process
and the coding excursion of the corresponding G-W tree converge in
law in the Skorokhod space of c\`adl\`ag paths. This work generalizes
Aldous' result \cite{Al} which concerns the Brownian case. In this
section we will prove that in the stable case, an invariance
principle also holds for sequences of G-W forests conditioned by
their mass.

Recall from section \ref{disc} that for an offspring distribution
$\mu$ we have set $\nu(k)=\mu(k+1)$, for $k=-1,0,1,\dots$.  We make
the following assumption:
\[\mbox{(H)}\;\;\;\;\left\{\begin{array}{l}
\mbox{$\mu$ is aperiodic and there is an increasing sequence
$(a_n)_{n\ge0}$}\\ \mbox{such that $a_n\rightarrow+\infty$ and
$S_n/a_n$ converges in law as $n\rightarrow+\infty$}\\\mbox{toward
the law of a non-degenerated r.v. $\theta$.}\end{array}\right.\]
Note that we are necessarily in the critical case, i.e. $\sum_k
k\mu(k)=1$, and that the law of $\theta$ is stable. Moreover, since
$\nu(-\infty,-1)=0$, the support of the L\'evy measure of $\theta$
is $[0,\infty)$ and its index $\alpha$ is such that $1<\alpha\le 2$.
Also $(a_n)$  is a regularly varying sequence with index $\alpha$.
Under hypothesis (H), it has been proved by Grimvall \cite{Gr} that
if $Z$ is the  G-W process associated to a tree or a forest with
offspring distribution $\mu$, then
\[\left(\frac1{a_n}Z_{[nt/a_n]},\,t\ge0\right)\Rightarrow(\overline{Z}_t,
\,t\ge0)\,,\;\;\;\mbox{as $n\rightarrow+\infty$,}\] where
$(\overline{Z}_t,\,t\ge0)$ is a continuous state branching process.
Here and in the sequel, $\Rightarrow$ will stand for the weak
convergence in the Skorokhod space of c\`adl\`ag trajectories.
Recall from section 2 the definition of the discrete process
$(S,H,C)$. Under the same hypothesis, Duquesne and Le Gall have
proved in \cite{DL}, Corollary 2.5.1 that
\begin{equation}\label{348}
\left[\left(\frac{1}{a_n}S_{[nt]},\frac{a_n}{n}H_{[nt]},
\frac{a_n}{n}C_{2nt}\right),\,t\ge0\right] \Rightarrow\left[(X_t,
\bar{H}_t,\bar{H}_t),\,t\ge0)\right]\,,\;\;\;\mbox{as
$n\rightarrow+\infty$,}\end{equation} where $X$ is a stable L\'evy
process with law $\theta$ and $\bar{H}$ is the associated height
process, as defined in section 3.

Again we fix a real $s>0$ and we consider a sequence of positive
integers $(k_n)$ such that
\begin{equation}\label{2104}
\frac{k_n}{a_n}\rightarrow s\,,\;\;\;\;\mbox{as
$n\rightarrow+\infty$.}\end{equation} Recall the notations of
section 2. For any $n\ge1$, let $(X^{br,n},\bar{H}^{br,n},C^{br,
n})$ be the process whose law is this of
\[\left[\left(\frac1{a_n}S_{[nt]},\frac{a_n}nH_{[nt]},
\frac{a_n}nC_{2nt}\right),\,0\le t\le1\right],\] under
$\mathbb{Q}_\mu(\,\cdot\,|\,\zeta(\tau_1)+\dots+\zeta(\tau_{k_n})=n)$.
Note that we could also define this three dimensional process over
the whole halfline $[0,\infty)$, rather than on $[0,1]$. However,
from the definitions in section 2, $\bar{H}^{br,n}$ and $C^{br, n}$
would simply vanish over $[1,\infty)$ and $X^{br,n}$ would be
constant. Here is the conditional version of the invariance
principle that we have recalled in (\ref{348}).

\begin{theorem}\label{main}
As $n$ tends to $+\infty$, we have
\begin{eqnarray*}
(X^{br,n},\bar{H}^{br,n},C^{br,n})\Longrightarrow
(X^{br},\bar{H}^{br},\bar{H}^{br})\,.
\end{eqnarray*}
\end{theorem}

\noindent In order to give a sense to the convergence of the L\'evy
forest, we may consider the trees ${\cal T}^{br,n}$ and  ${\cal
T}^{br}$ which are coded respectively  by the continuous processes
with compact support, $C^{br,n}_u$ and $\bar{H}^{br}_u$, in the
sense given at the beginning of section 3 (here we suppose that
these processes are defined on $[1,\infty)$ and both equal to 0 on
this interval). Roughly speaking the trees  ${\cal T}^{br,n}$ and
${\cal T}^{br}$ are obtained from the original (conditioned) forests
by rooting all the trees of these forests at a same root.
\begin{corollary}\label{cor}
The sequence of trees ${\cal T}^{br,n}$ converges weakly in the
space $\mathbb{T}$ endowed with the Gromov-Hausdorff topology
towards ${\cal T}^{br}$.
\end{corollary}
\noindent {\it Proof}. This results is a consequence of the weak
convergence of the contour function $C^{br,n}$ toward $\bar{H}^{br}$
and the inequality
\[d_{GH}({\cal T}_g,{\cal T}_{g'})\le 2\|g-g'\|\,,\]
which is proved in \cite{DL2}, see Lemma 2.3. (We recall that
$d_{GH}$ the Gromov-Hausdorff distance which has been defined in
section 3.)\QED

\noindent A first step for the proof of Theorem \ref{main} is to
obtain the weak convergence of $(X^{br,n},\bar{H}^{br,n})$
restricted to the Skorokhod space $\mathbb{D}([0,t])$ for any $t<1$.
Then we will derive the convergence on $\mathbb{D}([0,1])$ from an
argument of cyclic exchangeability. The convergence of the third
coordinate $C^{br,n}$ is a consequence of its particular expression
as a functional of the process $\bar{H}^{br,n}$. In the remainder of
the proof, we suppose that $S$ is defined on the same probability
space as $X$ and has step distribution $\nu$ under $\p$. Define also
$T_k=\inf\{i:S_i=-k\}$, for all integers $k\ge0$. Hence the process
$(X^{br,n},\bar{H}^{br,n},C^{br,n})$ has the same law as
\[\left[\left(\frac1{a_n}S_{[nt]},\frac{a_n}nH_{[nt]},
\frac{a_n}nC_{2nt}\right),\,0\le t\le1\right],\] under the
conditional probability $\p(\,\cdot\,|\,T_{k_n}=n)$.

\begin{lemma}\label{convweak}
For any $t<1$, as $n$ tends to $+\infty$, we have
\begin{eqnarray*}
\left[(X^{br,n}_u,\bar{H}^{br,n}_u),\,0\le u\le
t\right]\Longrightarrow \left[(X^{br}_u,\bar{H}^{br}_u),\,0\le y\le
t\right]\,.
\end{eqnarray*}
\end{lemma}
\noindent {\it Proof}.  From Feller's combinatorial lemma, see
\cite{Fe}, we have $\p(T_k=n)=\frac kn\p(S_n=-k)$, for all $n\ge1$,
$k\ge0$. Let $F$ be any bounded and continuous functional on
$\mathbb{D}([0,t])$. By the Markov property at time $[nt]$,
\begin{eqnarray}
&&\e[F(X^{br,n}_u,\bar{H}^{br,n}_u;0\le u\le
t)]=\e\left[F\left(\frac 1{a_n}S_{[nu]},\frac{a_n}nH_{[nu]};0\le
u\le
t\right)\,|\,T_{k_n}=n\right]\nonumber\\&=&\e\left(\ind_{\{[nt]\le
T_{k_n}\}}\frac{\p_{S_{[nt]}}(T_{k_n}=n-[nt])}{\p(T_{k_n}=n)} \times
F\left(\frac 1{a_n}S_{[nu]},\frac{a_n}nH_{[nu]};0\le u\le
t\right)\right)\nonumber\\
&=&\e\left(\ind_{\{\frac1{a_n}\underline{S}_{[nt]}\ge
-\frac{k_n}{a_n}\}}\frac{n(k_n+S_{[nt]})}{k_n(n-[nt])}\frac{\p_{S_{[nt]}}
(S_{n-[nt]}=-k_n)}{\p(S_n=-k_n)}\right.\label{density}\\&&\left.\qquad\qquad\qquad
\qquad\qquad\qquad\qquad\quad\times F\left(\frac
1{a_n}S_{[nu]},\frac{a_n}nH_{[nu]};0\le u\le
t\right)\right)\,.\nonumber
\end{eqnarray}
where $\underline{S}_k=\inf_{i\le k}S_i$. To simplify the
computations in the remainder of this proof, we set $P^{(n)}$ for
the law of the process $\left(\frac
1{a_n}S_{[nu]},\frac{a_n}nH_{[nu]};u\ge0\right)$ and $P$ will stand
for  the law of the process $(X_u,\bar{H}_u;u\ge0)$. Then
$Y=(Y^1,Y^2)$ is the canonical process of the coordinates on the
Skorokhod space $\mathbb{D}^2$ of c\`adl\`ag paths from $[0,\infty)$
into $\R^2$. We will also use special notations for the densities
introduced in (\ref{fpb}) and (\ref{density}):
\begin{eqnarray*}
D_t&=&\ind_{\{\underline{Y}^1_t\ge-s\}}\frac{s+Y^1_t}{s(1-t)}
\frac{\hat{p}_{1-t}(Y^1_t+s)}{\hat{p}_1(s)}\,,\;\;\;\mbox{and}\\
D_t^{(n)}&=&\ind_{\{\underline{Y}^1_{[nt]}\ge
-\frac{k_n}{a_n}\}}\frac{n(k_n+a_nY^1_{[nt]})}{k_n(n-[nt])}\frac{\p_{a_nY^1_{[nt]}}
(S_{n-[nt]}=-k_n)}{\p(S_n=-k_n)}\,,
\end{eqnarray*}
where $\underline{Y}^1_t=\inf_{u\le t}Y_u^1$. Put also $F_t$ for
$F(Y_u,0\le u\le t)$. To obtain our result, we have to prove that
\begin{equation}\label{5497}\lim_{n\rightarrow+\infty}
|E^{(n)}(F_tD^{(n)}_t)-E(F_tD_t)|=0\,.\end{equation} Let $M>0$ and
set $I_M(x)\eqdef \ind_{[-s,M]}(x)$. By writing
\[E^{(n)}(F_tD^{(n)}_t)=E^{(n)}(F_tD^{(n)}_tI_M(Y_t^1))+
E^{(n)}(F_tD^{(n)}_t(1-I_M(Y_t^1))\] and by doing the same for
$E(F_tD_t)$, we have the following upper bound for the term in
(\ref{5497})
\begin{eqnarray*}
|E^{(n)}(F_tD^{(n)}_t)-E(F_tD_t)|&\le&
|E^{(n)}(F_tD^{(n)}_tI_M(Y^1_t))-E(F_tD_tI_M(Y^1_t))|\\
&&+CE^{(n)}(D^{(n)}_t(1-I_M(Y^1_t)))+CE(D_t(1-I_M(Y^1_t)))\,,\end{eqnarray*}
where $C$ is an upper bound for the functional $F$. But since $D_t$
and $D_t^{(n)}$ are densities,  $E^{(n)}(D_t^{(n)})=1$ and
$E(D_t)=1$, hence
\begin{eqnarray}
|E^{(n)}(F_tD^{(n)}_t)-E(F_tD_t)|&\le&
|E^{(n)}(F_tD^{(n)}_tI_M(Y^1_t))-E(F_tD_tI_M(Y^1_t))|\label{label}\\
&&+C[1-E^{(n)}(D^{(n)}_tI_M(Y^1_t))]+C[1-E(D_tI_M(Y^1_t))]\,.\nonumber\end{eqnarray}
Now it remains to prove that the first term of the right hand side
of the inequality (\ref{label}) tends to 0, i.e.
\begin{equation}\label{remains}
|E^{(n)}(F_tD^{(n)}_tI_M(Y^1_t))-E(F_tD_tI_M(Y^1_t))|\rightarrow0\,,
\end{equation} as $n\rightarrow+\infty$. Indeed, suppose that
(\ref{remains}) holds, then by taking $F_t\equiv1$, we see that the
second term of the right hand side of (\ref{label}) converges
towards the third one. Moreover, $E(D_tI_M(Y^1_t))$ tends to 1 as
$M$ goes to $+\infty$. Therefore the second and the third terms in
(\ref{label}) tend to 0 as $n$ and $M$ go to $+\infty$.

Let us prove (\ref{remains}). From the triangle inequality and the
expression of the densities $D_t$ and $D^{(n)}_t$, we have
\begin{eqnarray}
&&|E^{(n)}(F_tD^{(n)}_tI_M(Y^1_t))-E(F_tD_tI_M(Y^1_t))|\le\sup_{x\in[-s,M]}
|g_n(x)-g(x)|+\nonumber\\
&&\qquad\qquad\qquad\qquad|E^{(n)}(F_tD_tI_M(Y^1_t))-E(F_tD_tI_M(Y^1_t))|\,,\label{last}
\end{eqnarray}
where $g_n(x)=\frac{n(k_n+x)}{k_n(n-[nt])}\frac{\p_{x}
(S_{n-[nt]}=-k_n)}{\p(S_n=-k_n)}$ and
$g(x)=\frac{s+x}{s(1-t)}\frac{p_{1-t}(x,-s)}{p_1(0,-s)}$. But thanks
to Gnedenko local limit theorem and the fact that
$k_n/a_n\rightarrow s$, we have
\[\lim_{n\rightarrow+\infty}\sup_{x\in[-s,M]}|g_n(x)-g(x)|=0\,.\]
Moreover, recall that from Corollary 2.5.1 of Duquesne and  Le Gall
\cite{DL},
\[P^{(n)}\Rightarrow P\,,\]
as $n\rightarrow+\infty$, where $\Rightarrow$ stands for the weak
convergence of measures on $\mathbb{D}^2$. Finally, note that the
discontinuity set of the functional $F_tD_tI_M(Y^1_t)$ is negligible
for the probability measure $P$ so that the last term in
(\ref{last}) tends to 0 as $n$ goes to $+\infty$. \QED

\noindent Then we will prove the tightness of the sequence,
$(X^{br,n},\bar{H}^{br,n})$. Define the height process associated to
any downward skip free chain $x=(x_0,x_1,\dots,x_n,\dots)$, i.e.
$x_0=0$ and $x_i-x_{i-1}\ge-1$, as follows:
\[H_n^{x}=\mbox{\rm card}\,\big\{i\in\{0,\dots,n-1\}:x_k=\inf_{i\le
j\le n}x_j\big\}\,.\] Let also $t(k)$ be the first passage time of
$x$ by $t(k)=\inf\{i:x_i=-k\}$ and for $n\ge k$, when $t(k)<\infty$,
define the shifted chain:
\[\theta_{t(k)}(x)_i=\left\{\begin{array}{ll}
x_{i+t(k)}+k,&\;\;\;\mbox{if $i\leq n-t(k)$}\\
x_{t(k)+i-n}+x_n+k,&\;\;\;\mbox{is $n-t(k)\leq i\leq
n$}\end{array}\right.,\;\;\;\;i=0,1,\dots,n\,,\] which consists in
inverting the pre-$t(k)$ and the post-$t(k)$ parts of $x$ and
sticking them together.
\begin{lemma}\label{det}
For any $k\ge0$,  we have almost surely
\[H^{\theta_{t(k)}(x)}=\theta_{t(k)}(H^{x})\,.\]
\end{lemma}
\noindent {\it Proof}. It is just a consequence of the fact that
$t(k)$ is a zero of $H^{x}$.\QED
\begin{lemma}\label{unif}
Let $u_{k_n}$ be a random variable which is uniformly distributed
over $\{0,1,\dots,k_n\}$ and independent of $S$. Under
$\p(\,\cdot\,|\,T(k_n)=n)$, the first passage time $T(u_{k_n})$ is
uniformly distributed over $\{0,1,\dots,n\}$.
\end{lemma}
\noindent {\it Proof}. It follows from elementary properties of
random walks that for all $k\in\{0,1,\dots,k_n\}$, under
$\p(\,\cdot\,|\,T(k_n)=n)$, the chain $\theta_{T(k_n)}(S)$ has the
same law as $(S_i,\,0\le i\le n)$. As a consequence, for all
$j\in\{0,1,\dots,n\}$,
\[P(T(k)=j\,|\,T(k_n)=n)=P(T(k_n-k)=n-j\,|\,T(k_n)=n)\,,\]
which allows us to conclude. \QED

\begin{lemma}\label{tightness}
The family of processes
\begin{eqnarray*}
(X^{br,n},\bar{H}^{br,n})\,,\;\;\;n\ge1
\end{eqnarray*}
is tight.
\end{lemma}
\noindent {\it Proof}. Let $\mathbb{D}([0,t])$ be the Skorokhod
space of c\`adl\`ag paths from $[0,t]$ to $\R$. In Lemma
\ref{convweak} we have proved the weak convergence of
$(X^{br,n},\bar{H}^{br,n})$ restricted to the space
$\mathbb{D}([0,t])$ for each $t>0$. Therefore, from Theorem 15.3 of
\cite{Bi}, it suffices to prove that for all $\delta\in(0,1)$ and
$\eta>0$,
\begin{equation}\label{23467}\lim_{\delta\rightarrow0}\limsup_{n\rightarrow+\infty}
\p\left(\sup_{s,t\in[1-\delta,1]}|X^{br,n}_{t}-X^{br,n}_{s}|>\eta,\,
\sup_{s,t\in[1-\delta,1]}|\bar{H}^{br,n}_{t}-
\bar{H}^{br,n}_{s}|>\eta\right)=0\,.\end{equation} Recall from Lemma
\ref{unif} the definition of the r.v. $u_{k_n}$ and define
$V_n=\inf\{t:X^{br,n}_t=-\frac kn\}$. Since from this lemma, $V_n$
is uniformly distributed over $\{0,1/n,\dots,1-1/n,1\}$, we have for
any $\varepsilon<1-\delta$,
\begin{eqnarray*}&&\p\left(\sup_{s,t\in[1-\delta,1]}|X^{br,n}_{t}-X^{br,n}_{s}|>\eta,\,
\sup_{s,t\in[1-\delta,1]}|\bar{H}^{br,n}_{t}-
\bar{H}^{br,n}_{s}|>\eta\right)\le\varepsilon+\delta+\\&&
\p\left(V_n\in[\varepsilon,1-\delta],\,
\sup_{s,t\in[1-\delta,1]}|X^{br,n}_{t}-X^{br,n}_{s}|>\eta,\,
\sup_{s,t\in[1-\delta,1]}|\bar{H}^{br,n}_{t}-
\bar{H}^{br,n}_{s}|>\eta\right)\,.\end{eqnarray*} Now for a
c\`adl\`ag  path $\omega$ defined on $[0,1]$ and $t\in[0,1]$, define
the shift:
\[\theta_{t}(\omega)_u=\left\{\begin{array}{ll}
\omega_{s+t}+u,&\;\;\;\mbox{if $s\leq 1-t$}\\
\omega_{t+u-1}+\omega_u+k,&\;\;\;\mbox{is $1-t\leq s\leq
1$}\end{array}\right.,\;\;\;\;u\in[0,1]\,,\] which consists in
inverting the paths $(\omega_u,\,0\le u\le t)$ and $(\omega_u,\,t\le
u\le 1)$ and sticking them together. We can check on a picture the
inclusion:
\begin{eqnarray*}
&&\{V_n\in[\varepsilon,1-\delta],\,
\sup_{s,t\in[1-\delta,1]}|X^{br,n}_{t}-X^{br,n}_{s}|>\eta,\,
\sup_{s,t\in[1-\delta,1]}|\bar{H}^{br,n}_{t}-
\bar{H}^{br,n}_{s}|>\eta\}\subset\\
&&\{
\sup_{s,t\in[0,1-\varepsilon]}|\theta_{V_n}(X^{br,n})_{t}-\theta_{V_n}(X^{br,n})_{s}|>\eta,\,
\sup_{s,t\in[0,1-\varepsilon]}|\theta_{V_n}(\bar{H}^{br,n})_{t}-
\theta_{V_n}(\bar{H}^{br,n})_{s}|>\eta\}\,.
\end{eqnarray*}
>From Lemma \ref{det} and the straightforward identity in law
$X^{br,n}\ed\theta_{V_n}(X^{br,n})$, we deduce the two dimensional
identity in law
$(X^{br,n},\Bar{H}^{br,n})\ed(\theta_{V_n}(X^{br,n}),\theta_{V_n}(\bar{H}^{br,n}))$.
Hence from the above inequality and inclusion,
\begin{eqnarray*}&&\p\left(\sup_{s,t\in[1-\delta,1]}|X^{br,n}_{t}-X^{br,n}_{s}|>\eta,\,
\sup_{s,t\in[1-\delta,1]}|\bar{H}^{br,n}_{t}-
\bar{H}^{br,n}_{s}|>\eta\right)\le\varepsilon+\delta+\\&& \p\left(
\sup_{s,t\in[0,1-\varepsilon]}|X^{br,n}_{t}-X^{br,n}_{s}|>\eta,\,
\sup_{s,t\in[0,1-\varepsilon]}|\bar{H}^{br,n}_{t}-
\bar{H}^{br,n}_{s}|>\eta\right)\,.\end{eqnarray*} But from Lemma
\ref{convweak} and Theorem 15.3 in \cite{Bi}, we have
\[\lim_{\delta\rightarrow0}\limsup_{n\rightarrow+\infty}
\p\left(\sup_{s,t\in[0,1-\varepsilon]}|X^{br,n}_{t}-X^{br,n}_{s}|>\eta,\,
\sup_{s,t\in[0,1-\varepsilon]}|\bar{H}^{br,n}_{t}-
\bar{H}^{br,n}_{s}|>\eta\right)=0\,.\] which yields
(\ref{23467}).\QED

\noindent {\it Proof of Theorem $\ref{main}$}. Lemma \ref{convweak}
shows that the sequence of processes $(X^{br,n},\bar{H}^{br,n})$
converges toward $(X^{br},\bar{H}^{br})$ in the sense of finite
dimensional distributions. Moreover tightness of this sequence has
been proved in Lemma \ref{tightness}, so we conclude from Theorem
15.1 of \cite{Bi}. The convergence of the two first coordinates in
Theorem \ref{main} is proved, i.e.
$(X^{br,n},\bar{H}^{br,n})\Longrightarrow (X^{br},\bar{H}^{br})$.
Then we may deduce the functional convergence of the third
coordinates from this convergence in law following similar arguments
as in Theorem 2.4.1 of \cite{DL} or in Theorem 3.1 of \cite{Du}:

>From (\ref{ch}), we can recover the contour process of $X^{br,n}$ as
follows set $b_i=2i-\Bar{H}^{br,n}_i$, for $0\leq i \le n$. For
$i\le n-1$ and $t\in[b_i,b_{i+1})$
\begin{equation*}
C^{br,n}_{t/2}=\left\{ \begin{array}{ll}
(\Bar{H}^{br,n}_i-(t-b_i))^+& \textrm{ if
}t\in[b_i,b_{i+1}-1),\\
(\Bar{H}^{br,n}_{i+1}-(b_{i+1}-t))^+, & \textrm{ if
}t\in[b_{i+1}-1,b_{i+1}),
\end{array}\right.\,.
\end{equation*}
Hence for $0\leq i\le n-1$,
\begin{equation}\label{descont}
\sup_{b_i\leq t<b_{i+1}}\left| C^{br,n}_{t/2}
-\Bar{H}^{br,n}_i\right|\leq \left| \Bar{H}^{br,n}_{i+1}
-\Bar{H}^{br,n}_i\right|+1\,.
\end{equation}
Now, we define $h_n(t)=i$, if $t\in[b_i, b_{i+1})$ and $i\le n-1$,
and $h_n(t)=n$, if $t\in[2n-2,2n]$. The definitions of $b_i$ and
$h_n$ implies
 \[
 \sup_{0\leq t\leq 2n}\left| h_n(t)-\frac{t}{2}\right|\leq \frac12
  \sup_{0\leq k\leq n}\Bar{H}^{br,n}_{k} +1.
 \]
Next, we set $f_n(t)=h_n(nt)/n$.  By (\ref{descont}), we have
\[
\sup_{0\leq t\leq 2}\frac{a_n}{n}\left| C^{br,n}_{nt/2}
-\Bar{H}^{br,n}_{nf_n(t)}\right|\leq \frac{a_n}{n}\sup_{0\leq t\leq
1}\left| \Bar{H}^{br,n}_{[nt]+1}
-\Bar{H}^{br,n}_{[nt]}\right|+\frac{a_n}{n},
\]
and
\[
\sup_{0\leq t\leq 2}\left| f_n(t)-\frac{t}{2}\right|\leq
\frac{1}{2a_n} \sup_{0\leq k\leq n}\frac{a_n}{n}\Bar{H}^{br,n}_{k}
+\frac{1}{p}.
\]
>From our hypothesis, we get
\[
\frac{a_n}{n}\sup_{0\leq t\leq 1}\left| \Bar{H}^{br,n}_{[nt]+1}
-\Bar{H}^{br,n}_{[nt]}\right|+\frac{a_n}{n}\to 0\quad \textrm{ as }
n\to \infty,
\]
and
\[
\frac{1}{2a_n} \sup_{0\leq k\leq n}\frac{a_n}{n}\Bar{H}^{br,n}_{k}
+\frac{1}{p}\to 0\quad \textrm{ as } n\to \infty,
\]
in probability. Hence, from the convergence of
$(X^{br,n},\bar{H}^{br,n})$ to $(X^{br},\bar{H}^{br})$, from Theorem
4.1 in \cite{Bi} and Skorokhod representation theorem we obtain the
convergence,
\[(X^{br,n},\bar{H}^{br,n},
C^{br,n}) \Rightarrow(X^{br},\bar{H}^{br},\bar{H}^{br}).
\]~\QED

\noindent{\bf Remarks}: By a classical time reversal argument, the
weak convergence of the first coordinate in Theorem \ref{main}
implies the main result of Bryn-Jones and R.A. Doney \cite{BD}.
Indeed, when $X$ is the standard Brownian motion, it is well known
that the returned first passage bridge $(s+X^{br}_{1-u},\,0\le u\le
1)$ is the bridge of a three dimensional Bessel process from 0 to
$s$ with length $1$. Similarly, the returned discrete first passage
bridge whose law is this of $(k_n+S_{n-i},\,0\le i\le n)$ under
$\p(\,\cdot\,|\,T(k_n)=n)$ has the same law as $(S_i,\,0\le i\le n)$
given $S_n=k_n$ and conditioned to stay positive. Then integrating
with respect to the terminal values and applying Theorem \ref{main}
gives the result contained in \cite{BD}.


\vspace*{.3in}

\begin{thebibliography}{99}

\bibitem{Al} \sc D.~Aldous: \rm
The continuum random tree. I. {\it Ann. Probab.}, {\bf 19}, No.1,
1-28, (1991).



\bibitem{Be} \sc J.~Bertoin (1996): {\it L\'evy Processes.} \rm
Cambridge University Press, Cambridge.

\bibitem{BCP} \sc J.~Bertoin, L.~Chaumont and J.~Pitman (2003): \rm Path
transformations of first passage bridges. \it Elect. Comm. in
Probab. \rm {\bf 8}, 155--166.

\bibitem{Bi} \sc P.~Billingsley: \rm Convergence of probability measures.
Second edition. \it John Wiley \& Sons, Inc., \rm New York, 1999.



\bibitem{BD} \sc  A.~Bryn-Jones and R.A.~Doney: \rm A functional
limit theorem for random walk conditioned to stay non-negative. {\it
J. London Math. Soc.} (2), {\bf 74}, no. 1, 244--258, (2006).

\bibitem{CU} \sc  L.~Chaumont and G.~Uribe: \rm Bridges of
self-similar Markov processes. {\it Work in progress}.

\bibitem{DMT} \sc A.~Dress, V.~Moulton and W.~Terhalle: \rm
$T$-theory: an overview. \it
 European J. Combin. \rm {\bf 17},  no. 2-3, 161--175, (1996).

\bibitem{Du} \sc T.~Duquesne: \rm
A limit theorem for the contour process of conditioned Galton-Watson
trees. {\it  Ann. Probab.}, {\bf 31}, No.2, 996-1027 (2003).

\bibitem{DL} \sc T.~Duquesne and J.F.~Le Gall: \rm
Random trees, L\'evy processes and spatial branching processes. {\it
Ast\'erisque $281$. Paris$:$ Soci\'et\'e Math\'ematique de France.
vi}, (2002).

\bibitem{DL3} \sc T.~Duquesne and J.F.~Le Gall: \rm
The Hausdorff measure of stable trees. {\it ALEA Lat. Am. J.
Probab.}, Math. Stat. {\bf 1}, 393--415, (2006).

\bibitem{DL2} \sc T.~Duquesne and J.F.~Le Gall: \rm
Probabilistic and fractal aspects of L\'evy trees. \it Probab.
Theory Related Fields, \rm {\bf 131}, no. 4, 553--603, (2005).

\bibitem{Ev} \sc S.~Evans: \rm Probability and trees.
{\it S\'eminaire de Probabilit\'es XXXVIII}, 1--4, \rm Lecture notes
in Math., 1857, Springer, Berlin, (2005).

\bibitem{EPW} \sc S.~Evans, J.~Pitman and A.~Winter: \rm
Rayleigh processes, real trees and root growth with re-grafting. \it
Probab. Theory Related Fields, \rm {\bf 134}, no 1, 81--126, (2006).

\bibitem{Fe}  \textit{\textrm{\textsc{W. Feller} An introduction to
probability theory and its applications. Vol. II. Second edition
John Wiley \& Sons, Inc., New York-London-Sydney, $1971$. }}

\bibitem{FPY} \sc P.J.~Fitzsimmons, J.~Pitman and M.~Yor (1993):
\rm Markovian bridges: construction, Palm interpretation, and
splicing. {\it Seminar on Stochastic Processes}, 1992, 101--134,
Progr. Probab., 33, Birkh\"auser Boston, MA.


\bibitem{Gr} \sc  A.~Grimvall: \rm
On the convergence of sequences of branching processes. {\it Ann.
Probab.}, {\bf 2}, 1027-1045 (1974).


\bibitem{Ji} \textsc{M.~Jirina}:
Stochastic branching processes with continuous state-space.
\textit{Czech. Math. J.}, \textbf{8}, 292-313, (1958).

\bibitem{Lab} \sc J.~Lamperti: \rm The limit of a sequence of branching processes. \it Z.
Wahrscheinlichkeitstheorie und Verw. Gebiete, \rm {\bf 7}, 271--288,
(1967).

\bibitem{Lab1} \sc J.~Lamperti: \rm Limiting distributions of branching processes. {\it
Proc. Fifth Berkeley Symp. Math. Statist. Probab.}, {\bf 2},
225--241. University California Press, Berkeley. (1967).

\bibitem{LG} \sc J.F.~Le Gall: \rm  Random real trees.
{\it Ann. Fac. Sci. Toulouse Math.}, (6), {\bf 15}, no. 1, 35--62,
(2006).

\bibitem{LL} \sc J.F.~Le Gall and Y.~Le Jan: \rm
Branching processes in L\'evy processes: The exploration process.
{\it Ann. Probab.}, {\bf 26}, No.1, 213--252, (1998).

\bibitem{Pi} \sc J.~Pitman: \rm Combinatorial stochastic processes.
\rm  Lectures from the 32nd Summer School on Probability Theory held
in Saint-Flour, July 7--24, 2002. Lecture Notes in Mathematics,
1875. {\it Springer-Verlag, Berlin}, 2006

\bibitem{WP} \sc J.~Pitman and M.~Winkel: \rm
Growth of the Brownian forest. {\it Ann. Probab.}, {\bf 33}, no. 6,
2188--2211, (2005).

\end{thebibliography}
\end{document}